\documentclass[11pt,a4paper,leqno]{article} 
\usepackage{a4wide}
\usepackage{latexsym}
\usepackage{amsmath}
\usepackage{amssymb}
\usepackage{stackrel}  
\usepackage{color}
\usepackage{graphicx}
\usepackage[linesnumbered,ruled,vlined]{algorithm2e}

\usepackage{hyperref}

\usepackage{authblk} 

\newtheorem{defin}{Definition}
\newtheorem{lemma}{Lemma}
\newtheorem{prop}{Proposition}
\newtheorem{theo}{Theorem}
\newtheorem{corol}{Corollary}
\newtheorem{example}{Example}
\newenvironment{proof}{\medskip\par\noindent{\bf Proof}}{\hfill $\Box$
\medskip\par}

\newcommand{\C}{\mathbb{C}}
\newcommand{\N}{\mathbb{N}}
\newcommand{\R}{\mathbb{R}}

\begin{document}
\title{On sequences preserving $q-$Gevrey asymptotic expansions}
\author[1]{Alberto Lastra}
\author[2]{S{\l}awomir Michalik}
\affil[1]{Universidad de Alcal\'a, Departamento de F\'isica y Matem\'aticas, Alcal\'a de Henares, Madrid, Spain. {\tt alberto.lastra@uah.es}}
\affil[2]{Faculty of Mathematics and Natural Sciences, College of Science. Cardinal Stefan Wyszy{\'n}ski University, Warszawa, Poland. {\tt s.michalik@uksw.edu.pl}}

\date{}

\maketitle
\thispagestyle{empty}
{ \small \begin{center}
{\bf Abstract}
\end{center}

The modification of the coefficients of formal power series is analyzed in order that such variation preserves $q-$Gevrey asymptotic properties, in particular $q-$Gevrey asymptotic expansions.  A characterization of such sequences is determined, providing a handy tool in practice. The sequence of $q-$factorials is proved to preserve $q-$Gevrey asymptotic expansions.

\bigskip

}
 \bigskip

\section{Introduction}

The present work is concerned with the sequences of positive real numbers preserving $q-$Gevrey asymptotic expansions. In~\cite{ichinobemichalik}, the authors noticed that the summability character of the formal solutions to certain $q-$difference-differential equations was preserved after the simplification of the $q-$difference operator ($0<q<1$). This is due to the property that the sequence $([p]_{q}!)_{p\ge0}$ preserves Gevrey asymptotic expansions. Roughly speaking, a sequence $m=(m_p)_{p\ge0}$ preserves Gevrey summability if for every formal power series $\hat{f}=\sum_{p\ge0}f_pz^p\in\C[[z]]$, $s>0$ and $d\in\R$ the series $\hat{f}$ is $s-$summable along direction $d\in\R$ if and only if $\mathcal{B}_{m}(\hat{f})=\sum_{p\ge0}\frac{f_p}{m_p}z^p$ is $s-$summable along $d$. In other words, the growth of the sequence $m$ is negligible with respect to the property of Gevrey summability of a formal power series. This concept coincides with that of preserving multisummability. In this respect, this paper is a continuation of~\cite{ichinobemichalik} to $q-$Gevrey asymptotic expansions.

The development of a parallel asymptotic theory in the $q-$Gevrey settings is more recent. We refer among others to the works and authors in~\cite{dreyfus,marottezhang,ramissauloyzhang,ramissauloyzhang2,tahara,viziozhang,zhang00}, and the references therein, and the application to the analysis of the asymptotic solutions to functional equations in~\cite{dreyfuseloy,dreyfuslastramalek,taharayamazawa0,lastramalek15,lastramaleksanz,ma22,tahara2,taharayamazawa0}, among others. Among the different $q-$analog approaches to the existing $q-$Gevrey asymptotic expansions and summability, we follow that developed in~\cite{zhang00}. In~\cite{viziozhang} the authors relate other approaches to the one considered in the present work (continuous and discrete definitions); see also~\cite{tahara} for a different approach in terms of $q-$exponential functions.

In this paper, we establish different kinds of spaces of formal power series which appear in asymptotic results associated to functional equations in the complex domain in a natural way. In~\cite{dreyfuslastramalek} $q-$difference-differential equations are studied. The analytic and formal solution admit decompositions into the sum of two terms which are associated to $q-$Gevrey asymptotic expansions of different orders (see Definition~\ref{defi0}). The appearance of a unique $q-$Gevrey asymptotic order can be seen, for example, in~\cite{ma17}. In~\cite{lama3}, two different asymptotic orders arise, one of them being Gevrey and the other $q-$Gevrey, in a sole singularly perturbed functional equation when exploring the asymptotic behavior of the analytic solution in the case that the values of the perturbation parameter become close to zero. Lastly, in~\cite{ma20} it is shown an asymptotic behavior of the analytic solutions to a $q-$difference$-$differential equation associated to mixed settings involving a $q-$Gevrey order governing a Gevrey suborder which stratifies the first one.

Regarding the variety of situations which may appear, one observes that the asymptotic expansions associated to $q-$difference equations are linked to a sole $q-$Gevrey order, or the coexistence of different $q-$Gevrey orders or even the existence of $q-$Gevrey orders modulated by Gevrey orders associated to one analytic solution in certain families of $q-$difference equations. On the other hand, the coexistence of different Gevrey asymptotic orders or associated to more general asymptotic expansions are being observed in~\cite{lamasa2}.

In the first part of the paper we analyze different spaces of formal power series whose coefficients are subject to $q-$Gevrey bounds, together with a sub-scale Gevrey. We prove that the sub-scale Gevrey is as precise as that linked to (in principal) more precise scales, linked to sequences satisfying properties naturally appearing in the asymptotic theory. Subsequently, we analyze the sequences of positive real numbers which preserve the $q-$Gevrey nature of a series of formal power series characterizing them as those of null $q-$Gevrey order, i.e. located between two Gevrey orders.

In~\cite{ichinobemichalik}, the concept of sequence preserving summability is put forward. The authors remark its interest by itself. It is proved there that all families of moment sequences of order 0 in the sense of Balser's moment summability~\cite{ba2} preserve summability. Indeed, both families coincide. The authors provide a characterization of such sequences, easier to check in practice. The present work gives a step forward in this theory defining the sequences that preserve q-Gevrey asymptotic expansions. The family of sequences preserving summability contains this new family. It is an open question whether this inclusion is strict or not. In this work, we have decided to adopt a $q-$analog definition of the previous notion in accordance with the nature of $q-$Gevrey asymptotic expansions. We say that a sequence of positive numbers $m$ preserves $q-$Gevrey asymptotic expansions if for every $s>0$, any $d\in\R$ and all formal power series $\hat{f}$ with coefficients in a complex Banach space, then the series $\hat{f}$ satisfies that its $q-$Borel transform of order $s$ converges on some neighborhood of the origin and its definition can be extended with adequate $q-$exponential growth along direction $d$ if and only if the same holds for the formal $m$-Borel transform of $\hat{f}$ (see Definition~\ref{defi7}). As a matter of fact, the concept of sequences preserving summability is equivalent to the parallel definition here adopted, due to Nevanlinna's Theorem (see~\cite{loday}, Theorem 5.3.9) which establishes that a series is summable for some positive Gevrey order $s$ along some direction $d$ if and only if its formal Borel transform of order $s$ converges on some neighborhood of the origin and such function can be extended to an infinite sector of bisecting direction $d$, with exponential growth at most $s$. In the $q-$Gevrey asymptotic framework we adopt, this definition is not equivalent to the existence of a unique function defined on some sector of bisecting direction $d$ and wide enough opening which admits the initial formal power series as its $q-$Gevrey asymptotic expansion  of order $s$, because unicity is not guaranteed. However, Laplace transform described in Theorem~\ref{teo164} guarantees that such $q-$Laplace transform of order $s$ admits the initial formal power series as its $q-$Gevrey asymptotic expansion of order $s$ on a sector of bisecting direction $d$ and opening larger that $2\pi$, i.e. on the Riemann surface of the logarithm.

For sequences preserving $q-$Gevrey asymptotic expansions, we establish an equivalence in a parallel way as that analyzed in~\cite{ichinobemichalik}, which is more easy to handle in practice (Theorem~\ref{teo2}). As a direct consequence of such characterization, one can affirm that the set of sequences preserving $q-$Gevrey asymptotic expansions is contained in that of sequences preserving summability in Corollary~\ref{coro1}. Finally, we prove that the sequence of $q-$factorials $([p]_{1/q}!)_{p\ge0}$ ($q>1$) does not only preserves summability as proved in~\cite{ichinobemichalik} but also preserves $q-$Gevrey asymptotic expansions.

The work is structured as follows: After recalling some known results and definitions on notions needed as a tool for our results such as some notions on $q-$calculus or $q-$Gevrey asymptotic expansions (Section~\ref{secpr}), we study some spaces of formal power series whose coefficients are upper bounded in terms of $q-$Gevrey and Gevrey sequences, determining properties related to them (Section~\ref{secsfs}). Section~\ref{sec3} is devoted to the study of spaces of sequences preserving $q-$properties. More precisely, we compare Gevrey sub-orders vs. other (in principal) more accurate sub-orders, and characterize the sequences preserving $q-$Gevrey order, as those sequences of $q-$Gevrey order 0 as those sequences between two Gevrey orders. The first main result (Theorem~\ref{teo2}) describes a characterization for a sequence of positive real numbers to preserve $q-$Gevrey asymptotic expansions easier to handle in practice. The work is concluded by proving that the sequence of $q-$factorials belongs to this family (Theorem~\ref{teopral}), which will be applied in practice to determine $q-$Gevrey asymptotic properties of the asymptotic solution to broad families of functional equations to be explored in the future.

\vspace{0.3cm}

\textbf{Notation:}

We write $\mathbb{N}=\{1,2,\ldots\}$ for the set of positive integers and $\mathbb{N}_0=\mathbb{N}\cup\{0\}$.

Throughout the present work, $\mathbb{E}$ stands for a complex Banach space with norm $\left\|\cdot\right\|_{\mathbb{E}}$. We write $\mathbb{E}[[t]]$ for the vector space of formal power series with coefficients in $\mathbb{E}$, and $\mathbb{E}\{t\}$ for the vector space of convergent power series on some neighborhood of the origin.

The open disc of positive radius $r$ is denoted by $D_r$. Given an open set $U\subseteq\C$, $\mathcal{O}(U,\mathbb{E})$ stands for the set of holomorphic functions defined in $U$ and values in $\mathbb{E}$. We simply write $\mathcal{O}(U)$ if $\mathbb{E}=\mathbb{C}$. For a sector $S$ of the complex domain with vertex at the origin, we write $T\prec S$ if $T$ is a subsector of $S$ with vertex at the origin such that $\overline{T}\setminus\{0\}\subseteq S\cap D_{\rho}$, for some $\rho>0$.

\section{Preliminary results}\label{secpr}

\subsection{Some notes on certain families of sequences and moment differentiation}

The asymptotic theory of formal solutions to functional equations rests on asymptotic expansion mainly related to Gevrey sequences, i.e. $(p!^{\alpha})_{p\ge0}$ for some fixed $\alpha>0$ (see~\cite{ba2,loday}). In recent years, a more general asymptotic theory has been developed involving sequences under certain appropriate assumptions providing the associated functional spaces with certain properties (see~\cite{sanzrev,jikalasa} and the references therein). Some of the properties satisfied by these sequences are the following:
\begin{itemize}
\item[(mg)] A sequence of positive real numbers $(M_p)_{p\ge0}$ is of moderate growth if there exists $A>0$ such that $M_{n+m}\le A^{n+m}M_nM_m$ for all $(n,m)\in\N_0^2$.
\item[(lc)] A sequence of positive real numbers $(M_p)_{p\ge0}$ is logarithmically convex if $M_{n}^2\le M_{n-1}M_{n+1}$ for every $n\ge1$.
\end{itemize}

We assume that the sequence is standardized by $M_0=1$.

In the present work, we show that the generality of such sequences would provide with different results than in the Gevrey case. However, the spaces of formal power series obtained coincide with those in which Gevrey sequences are considered (Proposition~\ref{prop1} and Proposition~\ref{prop4}).

We observe that the sequence involved in $q-$Gevrey asymptotic expansions $(q^{s\frac{p(p-1)}{2}})_{p\ge0}$ for some fixed $s>0$ is a sequence (lc), but not (mg). For such sequence, or the sequences appearing in Remark (a) page~\pageref{remark1}, different parallel asymptotic theories have arisen, as stated in the introduction, being that of Section~\ref{secsubsecqGev} the one considered in the present work.  

Given a sequence of positive real numbers $m=(m_p)_{p\ge0}$, we consider the formal operator
$$\partial_{m,z}\left(\sum_{p\ge0}\frac{a_p}{m_p}z^p\right)=\sum_{p\ge0}\frac{a_{p+1}}{m_p}z^p,$$
acting on $\mathbb{E}[[z]]$. These formal operators where initially described by W. Balser and M. Yoshino in~\cite{bayo}, and its definition can be naturally extended to holomorphic functions at some point. Recently, it has been proved to be extendable to sums and multisums of formal power series, as stated in~\cite{lamisu,lamisu3,lamisu4}. 

\subsection{Some notes on q-calculus}
We refer to~\cite{gara} for a reference on the following classical definitions and results, which are needed in the sequel.

Let $q>0$ such that $q\neq1$ be a real number. For every $\lambda\in\R$ the $q-$numbers are defined by $[\lambda]_{q}=\frac{q^{\lambda}-1}{q-1}$. Therefore, one has that $[n]_q=\frac{q^{n}-1}{q-1} $ for every $n\in\N$. The $q-$factorial is defined by $[0]_q!=1$ and for all $n\in\N$ one defines $[n]_{q}!=[1]_q[2]_q\cdots[n]_q$. It is well-known that $q-$factorials satisfy
\begin{equation}\label{e93}
[n]_q!=(1-q)^{-n}\prod_{p\ge0}\frac{1-q^{p+1}}{1-q^{p+1+n}},\qquad q\in(0,1),
\end{equation}
and 
$$[n]_{q}!=\prod_{p\ge0}\frac{1-q^{-1-p}}{1-q^{-n-1-p}}(q-1)^{-n}q^{\frac{(n+1)n}{2}},\qquad q>1.$$
Therefore, one has that
\begin{equation}\label{e90}
\lim_{n\to\infty}\frac{[n]_q!}{(1-q)^{-n}}=\prod_{p\ge0}1-q^{p+1}<\infty,\qquad q\in (0,1),
\end{equation}
and
\begin{equation}\label{e91}
\lim_{n\to\infty}\frac{[n]_q!}{q^{\frac{(n+1)n}{2}}(q-1)^{-n}}=\prod_{p\ge0}1-q^{-1-p}<\infty,\qquad q>1.
\end{equation}

It holds that
$$[p]_{q}!=[p]_{1/q}!q^{p(p-1)/2}\hbox{ for all }p\ge 0.$$

It will be useful to consider the notation of $q-$shift factorial for $q>1$ defined for every $n\in\N_0\cup\{\infty\}$ and $a\in\C$ by $(a;1/q)_{0}=1$, $(a;1/q)_{n}=\prod_{p=0}^{n-1}(1-\frac{a}{q^p})$ for $n\ge 1$ and $(a;1/q)_{\infty}=\prod_{p=0}^{\infty}(1-\frac{a}{q^p})<\infty$.

\subsection{Brief review on q-Gevrey asymptotic expansions}\label{secsubsecqGev}

In this section we recall the main facts about $q-$Gevrey asymptotic expansions of some positive order $s$. 

\begin{defin}\label{defi0}
Let $S$ be a sector with vertex at the origin. We say that $f\in\mathcal{O}(S,\mathbb{E})$ admits $\hat{f}(z)=\sum_{p\ge0}f_pz^p\in\mathbb{E}[[z]]$ as its $q-$Gevrey asymptotic expansion of order $s>0$ (at the origin) if for every $T\prec S$ there exist $C,A>0$ such that
$$\left\|f(z)-\sum_{p=0}^{N}f_pz^{p}\right\|_{\mathbb{E}}\le CA^Nq^{s\frac{N(N-1)}{2}}|z|^{N+1},$$
for all $N\in\N_0$ and all $z\in T$. 
\end{defin}

The $q-$Gevrey version of order $s>0$ of formal $q-$Borel transform is defined.

\begin{defin}
Let $\hat{f}(z)=\sum_{p\ge0}a_pz^p\in\mathbb{E}[[z]]$. The formal $q-$Borel transformation of order $s>0$ of $\hat{f}$ is defined by
$$\mathcal{B}_{q;s}(\hat{f})(z)=\sum_{p\ge0}\frac{a_p}{q^{s\frac{p(p-1)}{2}}}z^p.$$
\end{defin}

In order to assign $q-$sums to given formal power series, we now describe $q-$Laplace transform of order $k$. This procedure of $q-$summation of formal power series was firstly studied in~\cite{zhang00} for $s=1$, and has been extensively used in applications such as~\cite{lama00,ma17}, enumerated in the introduction of this work. We also refer to~\cite{tahara,viziozhang} for other operators in this respect. 

Jacobi Theta function is a holomorphic function on $\C^{\star}$ given by
$$\Theta(z)=\sum_{p\in\mathbb{Z}}\frac{1}{q^{\frac{p(p-1)}{2}}}z^p,$$
for every $z\in\C^{\star}$. For all $s>0$, we write $\Theta_{q^{s}}$ for the Jacobi Theta function associated with $q^{s}$ in place of $q$, it is to say $\Theta_{q^s}(z)=\sum_{p\in\mathbb{Z}}\frac{1}{q^{s\frac{p(p-1)}{2}}}z^p$. According to the factorization
$$\Theta_{q^s}(z)=\prod_{p\ge0}(1-\frac{1}{q^{s(p+1)}})(1+\frac{z}{q^{ps}})(1+\frac{1}{zq^{s(p+1)}}),\qquad z\in\C^{\star},$$
the set of zeroes of $\Theta_{q^s}$ is $\{-q^{ps}:p\in\mathbb{Z}\}$. The following is a slight modification of Lemma 4.1~\cite{lama00} which provides further information on the lower bounds of Jacobi Theta function.

\begin{lemma}[Lemma 4.1,~\cite{lama00}]\label{lema123}
Let $\delta>0$. There exists $C>0$ (which does not depend on $\delta$) such that 
$$|\Theta_{q^{s}}(z)|\ge C\delta\exp\left(\frac{\log^2|z|}{2s\log(q)}\right)|z|^{1/2},$$
for every $z\in\C^{\star}$ with $|1+z/q^{ks}|>\delta$ for every $k\in\mathbb{Z}$.
\end{lemma}

This growth property allows to define $q-$Laplace transform of order $1/s$ of functions satisfying $q-$exponential growth at infinity of order $1/s$ as follows.

\begin{defin}
Let $d\in\R$, $s>0$, and let $S$ be an infinite sector with vertex at the origin, bisecting direction $d$ and positive opening. Given $f\in\mathcal{O}(S,\mathbb{E})$, we say $f$ is of $q-$exponential growth of order $1/s$ if there exist $C,h>0$ and $\alpha\in\R$ such that
\begin{equation}\label{e131}
\left\|f(z)\right\|_{\mathbb{E}}\le C\exp\left(\frac{\log^2(|z|+h)}{2s\log(q)}\right)(|z|+h)^{\alpha},
\end{equation}
for every $z\in S$.
\end{defin}

\begin{defin}\label{def139}
Let $d\in\R$,$s>0$, $\rho>0$ and let $S$ be an infinite sector with vertex at the origin, bisecting direction $d$ and positive opening. We denote $\hat{S}=S\cup D_{\rho}$. We consider $f\in\mathcal{O}(\hat{S},\mathbb{E})$ of $q-$exponential growth of order $1/s$ and assume that (\ref{e131}) holds for $C,h>0,\alpha\in\R$ and $\rho$. We also assume that $f$ is continuous on $\overline{D_{\rho}}$. Let $\pi_{q^s}=\log(q)s$. The $q-$Laplace transform of order $1/s$ of $f$ along a direction $\gamma$ with $L_{\gamma}:=\mathbb{R}_+e^{i\gamma}\subseteq S$ is defined by
$$\mathcal{L}^{\gamma}_{q;s}(f)(z)=\frac{1}{\pi_{q^s}}\int_{L_{\gamma}}\frac{f(u)}{\Theta_{q^{s}}\left(\frac{u}{z}\right)}\frac{du}{u}.$$
\end{defin}

A deformation of the path $L_{\gamma}$ allows to substitute it by any other ray contained in $S$ in the previous definition.

The lower bounds showed in Lemma~\ref{lema123} guarantee analyticity of $q-$Laplace transform of order $1/s$ along suitable directions and applied on functions with $q-$exponential growth of order $1/s$.

\begin{prop}\label{prop157}
Let $\delta>0$ be as in Lemma~\ref{lema123}. We fix an infinite sector $S$ of bisecting direction $d\in\R$ and vertex at the origin. We denote $\hat{S}=S\cup D_{\rho}$ for some positive $\rho$. Let $f\in\mathcal{O}(\hat{S},\mathbb{E})$ be a function of $q-$exponential growth of order $1/s$ for some $s>0$ such that (\ref{e131}) holds for some $C,h>0$ and $\alpha\in\R$. Then, the function $\mathcal{L}_{q;s}^{\gamma}(f)$ given in Definition~\ref{def139} is a bounded holomorphic function on $\mathcal{R}_{\gamma,\delta}\cap D_{r}$ for any $0<r\le q^{-s(\alpha+1)}$, where
$$\mathcal{R}_{\gamma,\delta}=\left\{z\in\C^{\star}: \left|1+\frac{re^{i\gamma}}{z}\right|>\delta,\hbox{ for all }r\ge0\right\}.$$ 
\end{prop}

The next result in the case $s=1$ can be found in~\cite{zhang00}, and its proof can be adapted to the general framework. We have adapted the statement in our settings only providing the key steps for its proof for the sake of completeness.

\begin{theo}\label{teo164}
Let $\hat{f}\in\mathbb{E}[[z]]$ be such that $g=\mathcal{B}_{q;s}(\hat{f})\in\mathbb{E}\{z\}$ for some $s>0$. Assume that the function $g$ can be analytically extended into an infinite sector $S$ of bisecting direction $d\in\R$ and the extension is of $q-$exponential growth of order $1/s$. Then, for every $\gamma\in\R$ such that $L_{\gamma}=(0,\infty)e^{i\gamma}\subseteq S$ the function $\mathcal{L}_{q;s}^{\gamma}(g)$ admits $\hat{f}$ as its $q-$Gevrey asymptotic expansion of order $s$ on some finite sector of bisecting direction $d$ and opening $>2\pi$.
\end{theo}
\begin{proof}
Given $\hat{f}\in\mathbb{E}[[z]]$, we construct $g=\mathcal{B}_{q;s}(\hat{f})$. In virtue of Definition~\ref{def139} and Proposition~\ref{prop157}, the function $\mathcal{L}_{q;s}^{\gamma}(g)$ is a holomorphic and bounded function defined on a finite sector with bisecting direction $d$ and opening $>2\pi$ by analytically extending its definition with directions $\gamma$ which remain close to $d$.

Taking into account that
$$\frac{1}{\pi_{q^s}}\int_{L_{\gamma}}\frac{u^{p}}{\Theta_{q^s}\left(\frac{u}{z}\right)}\frac{du}{u}=z^{p}q^{s\frac{p(p-1)}{2}},$$
for every $z\in\mathcal{R}_{\gamma,\delta}$ and $p\in\N_0$ (see the proof of Lemma 3~\cite{ma20}) one can follow the classical proof in Gevrey asymptotic expansions slightly adapted to this framework (see for example Theorem 5.3.9~\cite{loday}) to conclude the expected asymptotic behavior of $\mathcal{L}_{q;s}^{\gamma}(g)$.
\end{proof}

\section{Spaces of formal power series}\label{secsfs}

\begin{defin}\label{defi1}
Let $q>0$ with $q\neq 1$, and fix $s\in\R$. We define the set $\mathbb{E}[[t]]_{q,s}$ consisting of all formal power series $\hat{u}(t)=\sum_{p\ge0}a_pt^p\in\mathbb{E}[[t]]$ such that there exist $A, B>0$, and $\alpha\in\R$ such that 
$$\left\|a_p\right\|_{\mathbb{E}}\le A B^{p}p!^{\alpha}q^{s\frac{p(p-1)}{2}},\qquad p\in\N_0.$$
\end{defin}

\noindent \textbf{Remark:}\label{remark1} 
\begin{itemize}
\item[(a)] The term $q^{s\frac{p(p-1)}{2}}$ can be substituted by $q^{s\frac{p^2}{2}}$. In view of the statement (\ref{e91}), it is straight to check that this term can also be substituted by $[p]_{q}!^s$ only in the case that $q>1$, due to the absence of symmetry observed in (\ref{e90}) with respect to (\ref{e91}).
\item[(b)] Observe the previous definitions can be reduced to $q>1$ as for $0<q<1$ the role of $s$ and $-s$ is interchanged.\label{remark}
\end{itemize}

The choice of a suborder Gevrey corresponds to practical reasons, due to their wide appearance in the summability of formal solutions to functional equations. This particular choice does not restrict the set in the following sense. Let us define $\mathbb{E}[[t]]_{q,s,\star}$ as the formal power series $\hat{u}(t)=\sum_{p\ge0}a_pt^p\in\mathbb{E}[[t]]$ such that there exist $A, B>0$, and a sequence $\mathbb{M}=(M_p)_{p\ge0}$ satisfying (mg) and (lc) properties, with
$$\left\|a_p\right\|_{\mathbb{E}}\le A B^{p}M_pq^{s\frac{p(p-1)}{2}},\qquad p\in\N_0,$$ 
then one has the following result.

\begin{prop}\label{prop1}
Let $q>0$ with $q\neq 1$, and fix $s\in\R$. Then, $\mathbb{E}[[t]]_{q,s}\equiv \mathbb{E}[[t]]_{q,s,\star}$
\end{prop}
\begin{proof}
Observe that the sequence $(p!^{\alpha})_{p\ge0}$ for $\alpha>0$ is a concrete sequence satisfying (mg) and (lc) conditions, so it is straight that $\mathbb{E}[[t]]_{q,s}\subseteq\mathbb{E}[[t]]_{q,s,\star}$. On the other hand, given any sequence $(M_p)_{p\ge0}$ which satisfies (mg) and (lc) properties is such that there exist $a_2,\delta>0$ with $M_p\le a_2^pp!^{\delta}$ for all $p\in\N$ (see Lemma 1.3.2.~\cite{thilliez}, for example). This entails that $\mathbb{E}[[t]]_{q,s,\star}\subseteq\mathbb{E}[[t]]_{q,s}$ holds.
\end{proof}

In contrast to Proposition~\ref{prop1}, there exist formal power series whose coefficients are linked to some $q-$Gevrey growth $s$ whose growth rate exceeds that of every element in $\mathbb{E}[[t]]_{q,s}$, as it is shown in the following example. 

\begin{example}
Let $q>1$ and $s>0$. The formal power series $\hat{u}(t)=\sum_{p\ge0}e^{p\log^2(p+1)}q^{s\frac{p(p-1)}{2}}t^p$ does not belong to $\mathbb{E}[[t]]_{q,s}$. Observe that the principal growth of the coefficients is related to the $q-$Gevrey sequence of order $s$. However, given any $\alpha\in\R$ and $A>0$ one has
$$\lim_{p\to\infty}\frac{e^{p\log^2(p+1)}q^{s\frac{p(p-1)}{2}}}{A^pp!^{\alpha}q^{s\frac{p(p-1)}{2}}}=+\infty.$$
The conclusion follows from here.
\end{example}

In an analogous way, we define the set of formal power series with coefficients in $\mathbb{E}$ subject to a fixed $q-$Gevrey order, and also a fixed Gevrey suborder as follows.

\begin{defin}\label{defi2}
Let $q>0$ with $q\neq 1$, and fix $s_1,s_2\in\R$. We define the set $\mathbb{E}[[t]]_{q,s_1}^{s_2}$ consisting of all formal power series $\hat{u}(t)=\sum_{p\ge0}a_pt^p\in\mathbb{E}[[t]]$ such that there exist $A, B>0$ with 
$$\left\|a_p\right\|_{\mathbb{E}}\le A B^{p}p!^{s_2}q^{s_1\frac{p(p-1)}{2}},\qquad p\in\N_0.$$
\end{defin}

The following properties hold regarding the previous sets of formal power series. They can be directly derived from the very definition of the sets in Definition~\ref{defi1} and Definition~\ref{defi2}. 

\begin{prop}\label{prop148}
Let $q>0$ with $q\neq 1$, and fix $s,s_1,s_2\in\R$. It holds that:
\begin{itemize}
\item[(i)] $\mathbb{E}[[t]]_{q,0}^{s}=\mathbb{E}[[t]]_s$, consisting of Gevrey sequences of order $s\in\R$. 

In particular, $\mathbb{E}[[t]]_{q,0}^{0}=\mathbb{E}\{t\}$. 
\item[(ii)] $$\bigcup_{h\in\R}\mathbb{E}[[t]]^{h}_{q,s_1}=\mathbb{E}[[t]]_{q,s_1}$$
\item[(iii)] $\mathbb{E}[[t]]_{q,s_2} $(resp. $\mathbb{E}[[t]]_{q,s_2}^{s_1}$) is a vector space, closed under formal differentiation. In the case that $\mathbb{E}$ is a Banach algebra and $s_2\ge0$, then $\mathbb{E}[[t]]_{q,s_2} $(resp. $\mathbb{E}[[t]]_{q,s_2}^{s_1}$) is a differential algebra.
\end{itemize} 
\end{prop}
\begin{proof}
All the properties follow from the very definition of the sets of formal power series. For the proof of (iii), one can make use of the classical estimates $1\le k!(p-k)!\le p!$ for every integer $0\le k\le p$, together with the fact that 
$$k(k+1)+(p-k)(p-k+1)=p(p-1)+k^2+(2-k)p,\qquad 0\le k\le p,$$
and $s_2\ge0$.
\end{proof}

\section{Spaces of sequences preserving q-properties}\label{sec3}

In this section, we consider new properties defining certain spaces of sequences of positive real numbers. It follows a similar structure as that of~\cite{ichinobemichalik}, where the authors describe and characterize sequences preserving Gevrey order and summability. In our framework, we enlarge the previous family allowing not only to maintain properties associated to Gevrey sequences, but also to preserve $q-$analogs of such properties. The different nature of $q-$Gevrey sequence causes crucial differences with respect to the previous approach.

In view of the second statement at the remark on page~\pageref{remark}, this previous assumption is made without loss of generality, treating the case $q\in (0,1)$ in the same way. For this reason, we will assume from now on that $q>1$ is a fixed real number. 

In addition to this, we will always assume that the sequence $m$ is normalized, in is to say, $m_0=1$.

We recall the definition of formal moment Borel operator, as described in~\cite{ba2}, Section 5.2. Although $m$ in the next definition can be any sequence of positive real numbers, we have maintained the word ``moment'' due to this sequence is usually a sequence of moments associated to some measure in practice.

\begin{defin}\label{defi7}
Let $m=(m_p)_{p\ge0}$ be a sequence of positive numbers. The formal $m$-Borel operator $\mathcal{B}_{m,t}:\mathbb{E}[[t]]\to\mathbb{E}[[t]]$ is defined by
$$\mathcal{B}_{m,t}\left(\sum_{p\ge0}a_pt^p\right)=\sum_{p\ge0}\frac{a_p}{m_p}t^p.$$
\end{defin}

\begin{defin}
A sequence $m=(m_{p})_{p\ge0}$ of positive real numbers is said to preserve $q-$Gevrey order if for every $s\in\R$ and all $\hat{u}\in\mathbb{E}[[t]]$, the following statements are equivalent:
\begin{itemize}
\item[(i)] $\hat{u}\in\mathbb{E}[[t]]_{q,s}$
\item[(ii)] $\mathcal{B}_{m,t}\hat{u}\in\mathbb{E}[[t]]_{q,s}$.
\end{itemize}
\end{defin}

In a similar way, one can define the sequences preserving $q-$Gevrey and Gevrey orders.

\begin{defin}\label{def176}
A sequence $m=(m_{p})_{p\ge0}$ of positive real numbers is said to preserve $q-$Gevrey and Gevrey orders if for every $s_1,s_2\in\R$ and all $\hat{u}\in\mathbb{E}[[t]]$, the following statements are equivalent:
\begin{itemize}
\item[(i)] $\hat{u}\in\mathbb{E}[[t]]_{q,s_1}^{s_2}$
\item[(ii)] $\mathcal{B}_{m,t}\hat{u}\in\mathbb{E}[[t]]_{q,s_1}^{s_2}$.
\end{itemize}
\end{defin}

Observe from (i) in Proposition~\ref{prop148} that the set of sequences of Definition~\ref{def176} is strictly contained in the set of sequences preserving Gevrey order, as given in~\cite{ichinobemichalik}, Definition 10, in which $s_1=0$. 

It is also worth mentioning that any sequence which preserves $q$-Gevrey and Gevrey orders is also a sequence which preserves $q-$Gevrey order. 

\begin{example}
Let $s\in\R$. The sequence $(p!^{s})_{p\ge0}$ preserves $q-$Gevrey order but it does not preserve $q-$Gevrey and Gevrey orders. The same holds for the sequence $(\Gamma(1+sp))_{p\ge0}$ for any fixed $s>0$.

Let $s\in\R$. The sequence $([p]_{1/q}!^s)_{p\ge0}$ preserves $q-$Gevrey and Gevrey orders. This is a direct consequence of the fact that $1\le [p]_{1/q}!\le \left(\frac{q}{q-1}\right)^{p}$ for every $p\in\N_0$.
\end{example}

\begin{defin}
A sequence $m=(m_p)_{p\ge0}$ is of $q-$Gevrey order $s\in\R$ if there exist $a,A>0$ and $\alpha,\beta\in\R$ such that
$$a^pp!^{\alpha}q^{s\frac{p(p-1)}{2}}\le m_{p}\le A^pp!^{\beta}q^{s\frac{p(p-1)}{2}},$$
for every $p\ge0$.
\end{defin}

As before, one can define the, in principle, wider class of sequences upper and lower bounded by a suborder determined by a more general sequence $\mathbb{M}=(M_p)_{p\ge0}$ which satisfies (mg) and (lc) properties. More precisely, one can state the following alternative definition.

\begin{defin}
A sequence $m=(m_p)_{p\ge0}$ is of (generalized) $q-$Gevrey order $s\in\R$ if there exist $a,A>0$ and a sequence of positive real numbers $\mathbb{M}=(M_p)_{p\ge0}$ satisfying (mg) and (lc) properties, such that
$$a^p\frac{1}{M_p}q^{s\frac{p(p-1)}{2}}\le m_{p}\le A^pM_pq^{s\frac{p(p-1)}{2}},$$
for every $p\ge0$.
\end{defin}

Nevertheless, the apparently generality is not really attained in view of the next result. 

\begin{prop}\label{prop4}
A sequence $m$ is of $q-$Gevrey order $s\in\R$ if and only if it is of (generalized) $q-$Gevrey order $s$.
\end{prop}
\begin{proof}
An analogous reasoning as that of the proof of Proposition~\ref{prop1} can be followed.
\end{proof}

At this point, one can provide with a characterization of sequences preserving $q-$Gevrey order. 

\begin{prop}
A sequence $m=(m_p)_{p\ge0}$ preserves $q-$Gevrey order if and only if $m$ is a sequence of $q-$Gevrey order 0, i.e. iff $m$ lies between two Gevrey sequences.
\end{prop}
\begin{proof}
Assume first that $m$ preserves $q-$Gevrey order, and consider $s\in\R$ together with the formal power series $\hat{u}(t)=\sum_{p\ge0}p!^sm_pt^p$. We observe that $\mathcal{B}_{m,t}\hat{u}\in\mathbb{C}[[t]]_{q,0}$. This entails that $\hat{u}\in\mathbb{C}[[t]]_{q,0}$. Therefore, there exist $\alpha\in\R$ and $A,B>0$ such that 
$$m_p \le AB^pp!^{\alpha-s},\qquad p\in\N_0.$$
Due to $m$ is a normalized sequence, it is straight to check that $m_p\le C^pp!^{\alpha-s}$ for every $p\ge0$, for some $C>0$. On the other hand, it holds that the formal power series $\hat{v}(t)=\sum_{p\ge0}p!^{s}t^p\in\C[[t]]_{q,0}$ which guarantees from the hypothesis that $\mathcal{B}_{m,t}\hat{v}\in\C[[t]]_{q,0}$. In other words, there exist $A,B>0$ and $\beta\in\R$ such that 
$$\frac{p!^{s}}{m_p}\le A B^{p}p!^{\beta},\qquad p\in\N_0,$$
or equivalently, $m_p\ge DE^pp!^{s-\beta}$, for some $D,E>0$, valid for all $p\in\N_0$. As the sequence $m$ is normalized, one can take $D=1$ for some $E>0$ large enough. The previous upper and lower bounds for $m_p$ for all $p\in\N_0$ guarantee that $m$ is of $q-$Gevrey order 0.

For the proof of the other implication, assume that $m=(m_p)_{p\ge0}$ is a $q-$Gevrey sequence of order 0, and choose $s\in\R$ and $\hat{u}(t)=\sum_{p\ge0}a_pt^p\in\mathbb{E}[[t]]$. First, assume that $\hat{u}\in\mathbb{E}_{q,s}[[t]]$. We will prove that $\mathcal{B}_{m,t}\hat{u}\in\mathbb{E}[[t]]_{q,s}$. From the hypothesis, it holds that $\left\|a_p\right\|_{\mathbb{E}}\le A_1B_1^pp!^{\alpha}q^{s\frac{p(p-1)}{2}}$ for some $A_1,B_1>0$ and some $\alpha\in\R$, valid for all $p\in\N_0$. As $m$ is a $q-$Gevrey sequence of order 0, it holds that
\begin{equation}\label{e240}
a^pp!^{\beta_1}\le m_p\le A^pp!^{\beta_2},\qquad p\in\N_0,
\end{equation}
for some $a,A>0$ and $\beta_1,\beta_2\in\R$.
Observe that for all $p\ge0$ one has
$$\frac{\left\|a_p\right\|_{\mathbb{E}}}{m_p}\le A_1\left(\frac{B_1}{a}\right)^pp!^{\alpha-\beta_1}q^{s\frac{p(p-1)}{2}}.$$
As a consequence, $\mathcal{B}_{m,t}\hat{u}\in\mathbb{E}[[t]]_{q,s}$. On the other hand, assume that $\mathcal{B}_{m,t}\hat{u}\in\mathbb{E}[[t]]_{q,s}$. Let us prove that $\hat{u}\in\mathbb{E}[[t]]_{q,s}$. From the hypothesis, one has the existence of $A_2,B_2>0$ and $\beta\in\R$ such that $\frac{\left\|a_p\right\|_{\mathbb{E}}}{m_p}\le A_2B_2^pp!^{\beta}q^{s\frac{p(p-1)}{2}}$ for every $p\in\N_0$. Together with (\ref{e240}) we conclude that 
$$\left\|a_p\right\|_{\mathbb{E}}=\frac{\left\|a_p\right\|_{\mathbb{E}}}{m_p}m_p\le A_2(AB_2)^pp!^{\beta+\beta_2}q^{s\frac{p(p-1)}{2}} $$
which yields the conclusion.
\end{proof}

As a matter of fact, the set of sequences preserving $q-$Gevrey and Gevrey orders coincides with that of sequences of null Gevrey order, it is to say, sequences $m=(m_p)_{p\ge0}$ such that $a^p\le m_p\le A^p$ for some $a,A>0$ for all $p\ge0$, in the terminology of Definition 2~\cite{ichinobemichalik}. The proof is skipped.

\begin{prop}
Let $m=(m_p)_{p\ge0}$ be a sequence of positive numbers. The following statements are equivalent:
\begin{itemize}
\item[(i)] The sequence $m$ preserves $q-$Gevrey and Gevrey orders
\item[(ii)] $m$ is a sequence of null Gevrey order.
\end{itemize}
\end{prop}

Observe that if a sequence $m=(m_p)_{p\ge0}$ preserves $q-$Gevrey and Gevrey orders, then it preserves Gevrey order in the sense of Definition 10~\cite{ichinobemichalik}, i.e. for every $s\in\R$ and any $\hat{u}\in\mathbb{E}[[t]]$ the equivalence $\hat{u}\in\mathbb{E}[[t]]_{q,0}^{s}$ if and only if $\mathcal{B}_{m,t}\hat{u}\in\mathbb{E}[[t]]_{q,0}^{s}$ holds.

The study of $q-$Gevrey asymptotic expansions is more involved than that of classical Borel summability. We consider different non-equivalent definitions in this direction. At this point, we restrict to the case of $q-$Gevrey asymptotic expansions of positive order: the case of $s=0$ is treated in~\cite{ichinobemichalik}, whereas the case $s<0$ deals with entire functions. 

\begin{defin}
A sequence $m$ is said to preserve $q-$Gevrey asymptotic expansions if for every $s>0$, $d\in\R$ and $\hat{u}\in\mathbb{E}[[t]]$ the following statements turn out to be equivalent:
\begin{itemize}
\item[(i)] $\mathcal{B}_{q;s}(\hat{u})\in\mathbb{E}\{t\}$, and this function can be extended on an infinite sector of bisecting direction $d$ with $q-$exponential growth of order $1/s$ on such sector.
\item[(ii)] $\mathcal{B}_{q;s}\mathcal{B}_{m,t}(\hat{u})\in\mathbb{E}\{t\}$, and this function can be extended on an infinite sector of bisecting direction $d$ with $q-$exponential growth of order $1/s$ on such sector.
\end{itemize}
\end{defin}

At this point, a characterization of sequences preserving $q-$Gevrey asymptotic expansions can be stated.

Observe that in order that given $\hat{u}\in\mathbb{E}[[t]]$ one has that $\mathcal{B}_{q;s}\hat{u}$ is convergent on some neighborhood of the origin, then $\hat{u}\in\mathbb{E}[[t]]_{q,s}^{0}$.

\begin{theo}\label{teo2}
A sequence $m=(m_p)_{p\ge0}$ preserves $q-$Gevrey asymptotic expansions if and only if for every $s>0$ and every $\theta\neq 0 \mod 2\pi$, $\mathcal{B}_{m,t}\left(\sum_{p\ge0}t^p\right)$ and $\mathcal{B}_{m^{-1},t}\left(\sum_{p\ge0}t^p\right)$ belong to $\mathbb{C}\{t\}$ and each of them can be extended to an infinite sector of bisecting direction $\theta$ with $q-$exponential growth of order $1/s$. 
\end{theo}
\begin{proof}
Assume that $m$ preserves $q-$Gevrey asymptotic expansions, and take $s>0$ and $\theta\neq 0 \mod 2\pi$. Consider the formal power series $\hat{u}(t)=\sum_{p\ge0}q^{s\frac{p(p-1)}{2}}t^p$. We observe that $\mathcal{B}_{q;s}\hat{u}(t)=\frac{1}{1-t}\in\mathbb{C}\{t\}$ and can be extended along every direction $\theta\neq0\mod 2\pi$ with $q-$exponential growth of order $1/s$. From the hypothesis, one has that $\mathcal{B}_{q;s}(\hat{u})\in\mathbb{C}\{t\}$ and it can be extended on an infinite sector of bisecting direction $d$ of $q-$exponential growth of order $1/s$ on such sector. From the hypothesis made, we have that the same holds for $\mathcal{B}_{q;s}\mathcal{B}_{m,t}(\hat{u})$ which coincides with 
$$\mathcal{B}_{q;s}\mathcal{B}_{m,t}(\hat{u})=\mathcal{B}_{m,t}\mathcal{B}_{q;s}(\hat{u})=\mathcal{B}_{m,t}\left(\sum_{p\ge0}t^p\right).$$
On the other hand, due to the hypothesis made, it also holds that $\mathcal{B}_{q;s}\mathcal{B}_{m^{-1},t}(\hat{u})\in\mathbb{C}\{t\}$ and this function can be extended on an infinite sector of bisecting direction $d$ with $q-$exponential growth of order $1/s$ on such sector. As one has that
$$\mathcal{B}_{q;s}\mathcal{B}_{m^{-1},t}(\hat{u})=\mathcal{B}_{m^{-1},t}\mathcal{B}_{q;s}(\hat{u})=\mathcal{B}_{m^{-1},t}\left(\sum_{p\ge0}t^p\right),$$
one can conclude the first part of the proof.

For the second part of the proof, we depart from $\hat{u}(t)=\sum_{p\ge0}u_{p}t^p\in\mathbb{E}[[t]]$ and assume that for every $s>0$ and $d\in\R$ one has that $\mathcal{B}_{q;s}(\hat{u})\in\mathbb{E}\{t\}$ and this function can be extended analytically to an infinite sector of bisecting direction $d$ and with $q-$exponential growth of order $1/s$ there. We will prove that the same holds for $\mathcal{B}_{q;s}\mathcal{B}_{m,t}(\hat{u})$. 

Let us consider the sequence $\tilde{m}=(m_pp!)_{p\ge0}$ and the Cauchy problem
$$(\partial_{\tilde{m},t}-\partial_z)\omega=0,$$
under the initial condition $\omega(0,z)=\mathcal{B}_{q;s}(\hat{u})(z)$. We observe that the formal solution of the previous problem is given by 
$$\hat{\omega}(t,z)=\sum_{p\ge0}\frac{(\mathcal{B}_{q;s}(\hat{u}))^{(p)}(z)}{m(p)p!}t^p.$$
Indeed, observe that
$$\mathcal{B}_{q;s}\mathcal{B}_{m,t}(\hat{u})=\sum_{p\ge0}\frac{u_p}{q^{s\frac{p(p-1)}{2}}m(p)}t^p=\sum_{p\ge0}\frac{(\mathcal{B}_{q;s}(\hat{u}))^{(p)}(0) }{p!m(p)}t^p=\hat{\omega}(t,0),$$
so $\hat{w}(t,0)$ defines a holomorphic function $\omega(t)$ defined on some neighborhood of the origin, say $D_r$ for some $r>0$. It rests to prove that $\hat{\omega}(t,0)$ can be prolonged to infinity with $q-$exponential growth along direction $d$. Let $0<r_0<r$. We have from Cauchy integral formula for the derivatives that
$$\omega(t,0)=\sum_{p\ge0}\frac{(\mathcal{B}_{q;s}(\hat{u}))^{(p)}(0)}{m_pp!}t^p=\frac{1}{2\pi i}\oint_{|\xi|=r_0}\frac{(\mathcal{B}_{q;s}(\hat{u}))(\xi)}{\xi}\psi(t/\xi)d\xi,$$
with $\psi(t)=\sum_{p\ge0}\frac{t^p}{m_p}=\mathcal{B}_{m,t}\left(\sum_{p\ge0}t^p\right)$. From the hypothesis, we have that $\psi(t)$ is convergent on some neighborhood of the origin, say $D_{r_1}$, and it can be extended analytically to an infinite sector of bisecting direction $d\neq0\mod 2\pi$ and with $q-$exponential growth of order $1/s$. We deform the integration path in
$$\omega(t)=\frac{1}{2\pi i}\oint_{|\xi|=r_0}\frac{(\mathcal{B}_{q;s}(\hat{u}))(\xi)}{\xi}\psi(t/\xi)d\xi$$
from $\{|\xi|=r_0\}$ to $\Gamma_1+\Gamma_2(R)$, where $\Gamma_1:(d+\epsilon,d+2\pi-\epsilon)\ni\theta\mapsto r_0e^{i\theta}$ for some small enough $\epsilon>0$ and $\Gamma_2(R)=\Gamma_{21}(R)+\Gamma_{22}(R)-\Gamma_{23}(R)$, with $\Gamma_{21}(R):(r_0,R)\ni s\mapsto se^{i(d+\epsilon)}$, $\Gamma_{23}(R):(r_0,R)\ni s\mapsto se^{i(d+2\pi-\epsilon)}$ 
and $\Gamma_{22}(R):(d-\epsilon,d+\epsilon)\ni \theta\mapsto Re^{i\theta}$. This yields holomorphy of $\omega(t)$ in $S_{d}\cup D_{r_0}$ by varying $R>0$.

Let $t\in S_{d}\cup D_{r_0}$. On the one hand, one has
\begin{multline*}\left\|\oint_{\Gamma_{1}}\frac{(\mathcal{B}_{q;s}(\hat{u}))(\xi)}{\xi}\psi(t/\xi)d\xi\right\|_{\mathbb{E}}\le\int_{d+\epsilon}^{d+2\pi-\epsilon}\left\|(\mathcal{B}_{q;s}(\hat{u}))(r_0e^{i\theta})\right\|_{\mathbb{E}}|\psi(t/(r_0e^{i\theta}))|d\theta \\
\le(2\pi-2\epsilon)\max_{|z|=r_0}\left\|\mathcal{B}_{q;s}(\hat{u}))(z)\right\|_{\mathbb{E}}C_1\exp\left(\frac{\log^2(|t|/r_0+\tilde{h})}{2s\log(q)}\right)\left(\frac{|t|}{r_0}+\tilde{h}\right)^{\alpha}
\end{multline*}
for some $C_1,\tilde{h}>0$ and $\alpha\in\R$. Observe that 
\begin{multline*}\log^{2}(|t|/r_0+\tilde{h})-\log^{2}(|t|+\tilde{h})=\log\left(\frac{|t|/r_0+\tilde{h}}{|t|+\tilde{h}}\right)\log((\frac{|t|}{r_0}+\tilde{h})(|t|+\tilde{h}))\\
\le \left[\max_{x>0}\log\left(\frac{x/r_0+\tilde{h}}{x+\tilde{h}}\right)\right]\log(\ell(|t|+\tilde{h})^2)\le \left[\max_{x>0}\log\left(\frac{x/r_0+\tilde{h}}{x+\tilde{h}}\right)\right](\log(\ell)+2\log(|t|+\tilde{h}))\\
\le C_3+C_4\log(|t|+\tilde{h}),
\end{multline*}
for some large enough integer $\ell>0$ and some $C_3,C_4>0$. We conclude that   
\begin{equation}\label{e393}\exp\left(\frac{\log^2(|t|/r_0+\tilde{h})}{2s\log(q)}\right)\le \exp\left(\frac{\log^2(|t|+\tilde{h})}{2s\log(q)}\right)\exp\left(\frac{C_3}{2s\log(q)}\right)(|t|+\tilde{h})^{\frac{C_4}{2s\log(q)}}.
\end{equation}
Observe moreover that 
\begin{equation}\label{e396}
\left(\frac{|t|}{r_0}+\tilde{h}\right)^{\alpha}\le C_5(|t|+\tilde{h})^{\alpha},
\end{equation}
for some $C_5>0$, valid for all $t$. From (\ref{e393}) and (\ref{e396}) one achieves that
\begin{equation}\label{e400}
\left\|\oint_{\Gamma_{1}}\frac{(\mathcal{B}_{q;s}(\hat{u}))(\xi)}{\xi}\psi(t/\xi)d\xi\right\|_{\mathbb{E}}\le \Delta_1\exp\left(\frac{\log^2(|t|+\tilde{h})}{2s\log(q)}\right)(|t|+\tilde{h})^{\alpha_1},
\end{equation}
with $\Delta_1=(2\pi-2\epsilon)\max_{|z|=r_0}\left\|\mathcal{B}_{q;s}(\hat{u}))(z)\right\|_{\mathbb{E}}C_1\exp\left(\frac{C_3}{2s\log(q)}\right)C_5$ and $\alpha_1=\alpha+\frac{C_4}{2s\log(q)}$.

Secondly, analogous estimates yield 
$$\left\|\oint_{\Gamma_{21}(R)}\frac{(\mathcal{B}_{q;s}(\hat{u}))(\xi)}{\xi}\psi(t/\xi)d\xi\right\|_{\mathbb{E}}\le \int_{r_0}^{R}\frac{1}{h}\left\|(\mathcal{B}_{q;s}(\hat{u}))(he^{i(d+\epsilon)})\right\|_{\mathbb{E}} |\psi(t/(he^{i(d+\epsilon)}))|dh.$$

Let $R=\frac{2|t|}{r_0}$. Taking into account the $q-$exponential growth of the terms in the integrand and the fact that for every $h\in[r_0,R]$ and large enough $R$ one has that 
\begin{multline*}
\log^2(h)+(\log(|t|/h))^{2}=2\log^2(h)+\log^2|t|-2\log|t|\log(h)=\log^2|t|+2\log(h)\log(h/|t|)\\
\le\log^2|t|+2\log(2|t|/r_0)\log(2/r_0),
\end{multline*}
one concludes that
\begin{equation}\label{e413}
\left\|\oint_{\Gamma_{21}(R)}\frac{(\mathcal{B}_{q;s}(\hat{u}))(\xi)}{\xi}\psi(t/\xi)d\xi\right\|_{\mathbb{E}}\le\Delta_2\exp\left(\frac{\log^2(|t|+\tilde{\delta}_2)}{2s\log(q)}\right)(|t|+\tilde{h}_2)^{\alpha_2}, 
\end{equation}
for some $\Delta_2,\tilde{\delta}_2,\tilde{h}_2>0$, and some $\alpha_2\in\R$.
The upper estimates for the integral along $\Gamma_{23}$ are upper bounded in the same manner, arriving at
\begin{equation}\label{e414}
\left\|\oint_{\Gamma_{23}(R)}\frac{(\mathcal{B}_{q;s}(\hat{u}))(\xi)}{\xi}\psi(t/\xi)d\xi\right\|_{\mathbb{E}}\le\Delta_3\exp\left(\frac{\log^2(|t|+\tilde{\delta}_3)}{2s\log(q)}\right)(|t|+\tilde{h}_3)^{\alpha_3}, 
\end{equation}
for some $\Delta_3,\tilde{\delta}_3,\tilde{h}_3>0$, and some $\alpha_3\in\R$. Finally, observe following analogous arguments that
\begin{multline*}
\left\|\oint_{\Gamma_{22}(R)}\frac{(\mathcal{B}_{q;s}(\hat{u}))(\xi)}{\xi}\psi(t/\xi)d\xi\right\|_{\mathbb{E}}\le \frac{\epsilon R}{\pi R} C_6(R+\tilde{\delta}_{4})^{\alpha_4}\exp\left(\frac{\log^2(R+\tilde{\delta}_4)}{2s\log(q)}\right)C_7  (\frac{|t|}{R}+\tilde{\delta}_{5})^{\alpha_5}\\
\times\exp\left(\frac{\log^2(\frac{|t|}{R}+\tilde{\delta}_5)}{2s\log(q)}\right)\le \Delta_4(\frac{2|t|}{r_0}+\tilde{\delta}_{4})^{\alpha_4}\exp\left(\frac{\log^2(\frac{2|t|}{r_0}+\tilde{\delta}_4)}{2s\log(q)}\right), 
\end{multline*}
for some $C_6,C_7,\tilde{\delta}_4,\tilde{\delta}_5>0$ and $\alpha_4,\alpha_5\in\R$, and where
$$\Delta_4=\epsilon C_6 C_7  (\frac{r_0}{2}+\tilde{\delta}_{5})^{\alpha_5}
\exp\left(\frac{\log^2(\frac{r_0}{2}+\tilde{\delta}_5)}{2s\log(q)}\right).$$
Analogous estimates as above allow us to conclude from the previous upper estimate together with (\ref{e400}), (\ref{e413}) and (\ref{e414}). 

The proof in the case of $m^{-1}$ can be developed in an analogous manner.
\end{proof}

As a matter of fact, any sequence preserving $q-$Gevrey asymptotic expansions preserves summability in the sense of~\cite{ichinobemichalik}.

\begin{defin}
Let $k>0$ and $d\in\R$. A formal power series $\hat{u}(t)=\sum_{p\ge0}a_pt^p\in\mathbb{E}[[t]]$ is $k-$summable along direction $d$ if $\mathcal{B}_{\Gamma_{1/k},t}\hat{u}=\sum_{p\ge0}\frac{a_p}{\Gamma(1+\frac{p}{k})}t^p$ defines a convergent function on some neighborhood of the origin which can be analytically prolonged as an analytic function, say $v$, into an infinite sector of bisecting direction $d$, say $S$, and such extension is of $k$ exponential growth, i.e. there exist $A,B>0$ such that $\left\|v(t)\right\|_{\mathbb{E}}\le Ae^{B|t|^{k}}$ for all $t\in S$.
\end{defin}

\begin{defin}[Definition 11,~\cite{ichinobemichalik}]\label{teo3}
A sequence of positive real numbers $m$ preserves summability if for every $k>0$, $d\in\R$ and every $\hat{u}\in\mathbb{E}[[t]]$, the following statements are equivalent:
\begin{itemize}
\item[(i)] $\hat{u}\in\mathbb{E}\{t\}_{k,d}$.
\item[(ii)] $\mathcal{B}_{m,t}\hat{u}\in\mathbb{E}\{t\}_{k,d}$.
\end{itemize}
\end{defin}

\begin{theo}[Theorem 1,~\cite{ichinobemichalik}]\label{teor4}
A sequence of positive real numbers $m$ preserves summability if and only if for every $k>0$ and every $\theta\neq 0 \mod 2\pi$, $\mathcal{B}_{m,t}\left(\sum_{p\ge0}t^p\right)$ and $\mathcal{B}_{m^{-1},t}\left(\sum_{p\ge0}t^p\right)$ belong to $\mathbb{C}\{t\}$ and they can be extended to an infinite sector of bisecting direction $\theta$ with exponential growth of order $k$. 
\end{theo}

 In view of Theorem~\ref{teo2}, one has the following result.

\begin{corol}\label{coro1}
A sequence of positive real numbers $m$ which preserves $q-$Gevrey asymptotic expansions is a sequence preserving summability.
\end{corol}
\begin{proof}
Given $m$ preserving $q-$Gevrey asymptotic expansions, it holds from Theorem~\ref{teo2} that for $s>0$ and $\theta\neq 0 \mod 2\pi$, the formal power series $\mathcal{B}_{m,t}\left(\sum_{p\ge0}t^p\right)$ and $\mathcal{B}_{m^{-1},t}\left(\sum_{p\ge0}t^p\right)$ are convergent on some neighborhood of the origin, and both can be extended to an infinite sector of bisecting direction $\theta$, say $S_{\theta}$, with $q-$exponential growth of order $s$. We observe that a function with $q-$exponential growth of order $s$ on $S_{\theta}$ satisfies bounds as in (\ref{e131}), and therefore is of $k-$exponential growth on $S_{\theta}$ for every $k>0$. This entails from Theorem~\ref{teor4} that $m$ preserves summability.
\end{proof}

\begin{example}
Given $A>0$, the sequence $(A^{p})_{p\ge0}$ preserves summability and $q-$Gevrey asymptotic expansions.
\end{example}

\begin{example}
The sequence $m=\left(\frac{(2p)!}{p!^2}\right)_{p\ge0}$ preserves $q-$Gevrey asymptotic expansions. Indeed, take $d\in\R$ with $d\neq 0\mod 2\pi$ and $s>0$. It holds that
$$\mathcal{B}_{m^{-1},t}\left(\sum_{p\ge0}t^p\right)=\sum_{p\ge0}\frac{(2p)!}{p!^2}t^p=\frac{1}{(1-4t)^{1/2}},$$
together with
$$\mathcal{B}_{m,t}\left(\sum_{p\ge0}t^p\right)=\sum_{p\ge0}\frac{p!^2}{(2p)!}t^p=\frac{1}{1-\frac{t}{4}}\left(1+\frac{1}{2}\left(\frac{t}{1-\frac{t}{4}}\right)^{1/2}\arcsin(\frac{1}{2}t^{1/2})\right).$$
Both functions have $q-$exponential growth along direction $d$. The conclusion follows.
\end{example}

\begin{theo} \label{teopral}
The sequence $([p]_{1/q}!)_{p\ge0}$ preserves summability and it also preserves $q-$Gevrey asymptotic expansions. 
\end{theo}
\begin{proof}
The fact that the sequence $m_{q}:=([p]_{1/q}!)_{p\ge0}$ preserves summability is stated and proved in Theorem 2~\cite{ichinobemichalik}. We now give proof to the fact that $([p]_{1/q}!)_{p\ge0}$ preserves $q-$Gevrey asymptotic expansions. We prove the characterization stated in Theorem~\ref{teo2}. Let $s>0$ and $\theta\neq 0\mod 2\pi$. We consider the formal power series
$$\hat{x}_1(t)=\sum_{p\ge0}\frac{1}{[p]_{1/q}!}t^p,\qquad \hat{x}_2(t)=\sum_{p\ge0}[p]_{1/q}!t^p.$$
We will prove that both formal power series are convergent on some neighborhood of the origin and admit analytic continuation on an infinite sector $S$ with bisecting direction $\theta$ and $q-$exponential growth of order $s$ on such sector. 

We observe that $\hat{x}_1$ coincides with the $q-$exponential function $e_{1/q}$. It is well-known (see Lemma 2.6~\cite{zhang01} or Proposition 2.2~\cite{tahara}) that $e_{1/q}$ is holomorphic on $D_{q/(q-1)}$ and can be analytically prolonged to $\C\setminus\R_{+}$ in such a way that for every sector $S_1$ with bisecting direction $\pi$ and opening smaller than $\pi$ it holds that
$$|e_{1/q}(z)|\le K,$$
for some $K>0$ (which depends on $S_1$) and any $z\in S_1$ with $|z|\ge1$. Indeed, $K_1$ can be substituted by a $q-$exponential decreasement at infinity. This allows to conclude the proof for $\hat{x}_1$. 

We now consider the formal power series $\hat{x}_2$. We consider the initial value problem
\begin{equation}\label{e487}
\left\{ \begin{array}{l}
             (\partial_{\tilde{m},t}-\partial_{z})u=0 \\
             u(t,0)=\frac{1}{1-t}
             \end{array}
   \right.,
\end{equation}
where $\tilde{m}=(p![p]_{1/q}!)_{p\ge0}$. We observe that
$$\hat{u}(t,z)=\sum_{p\ge0}\frac{\partial_{\tilde{m},t}^{p}\left(\frac{1}{1-t}\right)}{p!}z^{p}$$ 
is the formal solution to (\ref{e487}). 
Since
$$1\le [p]_{1/q}!\le\left(\frac{q}{q-1}\right)^{p},\qquad p\ge0,$$
by Proposition  1~\cite{lamisu2}, there exists $r,A,B>0$ such that  
$$\sup_{|t|<r}\left|\partial_{\tilde{m},t}^{p}\left(\frac{1}{1-t}\right)\right|\le A B^p p!,\qquad p\ge0.$$
Hence 
$\hat{u}$ converges on the neighborhood of the origin $D_r\times D_{1/B}\subseteq\C^2$. We write $u(t,z)$ for the analytic function that this series defines. As a consequence of moment derivation applied on holomorphic functions at a point of holomorphy (i.e. $\partial_{m,t}^{p}\psi(0)/m_p=\psi^{(p)}(0)/p!$ for all $p\ge0$, any moment sequence and all $\psi\in\mathbb{E}\{t\}$), we observe that for small enough $t$, one has
$$u(0,z)=\sum_{p\ge0}\frac{\partial_{\tilde{m},t}^{p}\left(\frac{1}{1-t}\right)|_{t=0}}{p!}z^{p}=\sum_{p\ge0}\frac{\tilde{m}_p\left(\frac{1}{1-t}\right)^{(p)}|_{t=0}}{p!^2}z^{p}=\sum_{p\ge0}[p]_{1/q}!z^p.$$
Therefore, $\hat{x}_2$ defines a holomorphic function on some neighborhood of the origin. Moreover, taking into account (\ref{e93}), and from Cauchy integral formula one can write
$$u(0,z)=\frac{1}{2\pi i}\oint_{|\nu|=\rho}\frac{1}{\nu(1-\nu)}\sum_{n\ge0}(1/q;1/q)_{n}\left(\frac{z}{(1-1/q)\nu}\right)^{n}d\nu,$$
for small enough $\rho>0$, valid on some neighborhood of $z=0$. At this point, similarly to Lemma 3~\cite{ichinobemichalik}, we apply Heine's transformation formula (\cite{gara}, Section 1.4; see also Proposition 2~\cite{ichinobemichalik}) to write the previous expression in the form 
$$\frac{1}{2\pi i}\oint_{|\nu|=\rho}\frac{1}{\nu(1-\nu)}\frac{(1/q;1/q)_{\infty}(z/((q-1)\nu);1/q)_{\infty}}{(qz/((q-1)\nu);1/q)_{\infty}}\sum_{p\ge0}\frac{(qz/((q-1)\nu);1/q)_{p}}{(z/((q-1)\nu);1/q)_{p}(1/q;1/q)_{p}}\frac{1}{q^p}d\nu.$$
For every $z\neq0$, the function $\nu\mapsto \frac{1}{(qz/((q-1)\nu);1/q)_{\infty}}$ is meromorphic in $\C$ with simple poles at $\nu_p=\frac{qz}{(q-1)q^p}$ for every $p\ge0$. The residue principle yields 
\begin{multline*}
u(0,z)=(1/q;1/q)_{\infty}\sum_{p\ge0}\frac{1}{1-\frac{qz}{(q-1)q^p}}\hbox{Res}_{\nu=\nu_p}\frac{1}{(qz/((q-1)\nu);1/q)_{\infty}}\frac{1}{\nu}\\
\hfill\times(q^{p-1};1/q)_{\infty}\sum_{j\ge0}\frac{(q^p;1/q)_{j}}{(q^{p-1};1/q)_{j}(1/q;1/q)_{j}q^j}\\
=(1/q;1/q)_{\infty}\sum_{p\ge0}\frac{1}{1-\frac{qz}{(q-1)q^p}}\hbox{Res}_{\nu=\nu_p}\frac{1}{(qz/((q-1)\nu);1/q)_{\infty}}\frac{1}{\nu}\frac{(q^{p-1};1/q)_{\infty}(q^p;1/q)_{p}}{(q^{p-1};1/q)_{p}(1/q;1/q)_{p}q^p}\\
=(1/q;1/q)_{\infty}\sum_{p\ge0}\frac{1}{1-\frac{qz}{(q-1)q^p}}\hbox{Res}_{\nu=\nu_p}\frac{1}{(qz/((q-1)\nu);1/q)_{\infty}}\frac{1}{\nu}\frac{(1/q;1/q)_{\infty}(q^p;1/q)_{p}}{(1/q;1/q)_{p}q^p}\\
=(1/q;1/q)_{\infty}\sum_{p\ge0}\frac{1}{1-\frac{qz}{(q-1)q^p}}\frac{1}{(q^p;1/q)_p(1/q;1/q)_{\infty}}\frac{(1/q;1/q)_{\infty}(q^p;1/q)_{p}}{(1/q;1/q)_{p}q^p},
\end{multline*}
from which we conclude that
$$u(0,z)=\sum_{p\ge0}\frac{1}{1-\frac{qz}{(q-1)q^p}}\frac{(1/q;1/q)_{\infty}}{(1/q;1/q)_{p}}\frac{1}{q^p}.$$
Observe that $(1/q;1/q)_{\infty}\le (1/q;1/q)_{p}$ for all $p\ge0$, $\sum_{p\ge0}1/q^p<\infty$, and
$$\left|1-\frac{qz}{(q-1)q^p}\right|\ge 1-\frac{q|z|}{q-1}=1-\frac{|z|}{1-1/q}.$$
The function $u(0,z)$ can be analytically prolonged in $\C\setminus\R_{+}$ and it has polynomial growth at infinity (therefore $q-$exponential growth of order $s$) along direction $\theta\neq0\mod 2\pi$.
\end{proof}

\end{document}